\documentclass[12pt]{article}

\usepackage{amsmath, amsfonts, amssymb}

\usepackage[english]{babel}
\usepackage[utf8]{inputenc}
\usepackage{todonotes}
\usepackage{comment}
\usepackage[left=2.5cm, right=3cm, top=2cm, bottom=3cm, bindingoffset=0cm, marginparwidth=3cm]{geometry}

\usepackage{hyperref}
\usepackage{enumitem}
\usepackage{amsthm}

\usepackage{tikz}
\usetikzlibrary{arrows,calc}
\tikzset{
>=stealth',
}

\newcounter{theorems}
\newtheorem{prop}[theorems]{Proposition}
\newtheorem{thm}[theorems]{Theorem}

\newtheorem{mytheorem}{Theorem}

\theoremstyle{definition}
\newtheorem{defin}[theorems]{Definition}
\theoremstyle{remark}

\newtheorem{rem}[theorems]{Remark}
\newtheorem*{rem*}{Remark}

\makeatletter
\def\blfootnote{\gdef\@thefnmark{}\@footnotetext}
\makeatother

\def\Om {\Omega }

\def\e {\varepsilon }

\def\phi {\varphi}

\def\be {\begin{equation}}

\def\ee {\end{equation}}

\def\bt {\begin{thm}}

\def\et {\end{thm}}

\def\N {\mathcal N}

\def\F {\mathcal F}

\def\RR {\mathbb R}

\def\ZZ {\mathbb Z}

\def\nn {\mathbb N}

\def\QQ {\mathbb Q}

\def\orc{\mathcal L_+}

\def\TT{{\mathbb{T}^2}}

\def\Si{\Sigma}

\def\U{\mathcal U}

\def\Bif{\mathrm{Bif}}

\def\Si{\Sigma}

\def\hSi{{\hat{\Sigma}}}

\def\t{{\theta}}

\def\ddt{\frac{\partial}{\partial \t}}

\begin{document}

\title{Sparkling saddle loops of vector fields on surfaces}
\author{Ivan Shilin\footnote{National Research University Higher School of Economics, Moscow}}
\date{}

\maketitle
\begin{abstract}
We study bifurcations of vector fields on 2-manifolds with handles in generic one-parameter families unfolding vector fields with a separatrix loop of a hyperbolic saddle. These bifurcations can differ drastically from the analogous bifurcations on the sphere. The reason is that, on a surface, a free separatrix of a hyperbolic saddle may wind toward the separatrix loop of the same saddle. When this loop is broken, sparkling saddle loops emerge. In the orientable case, the parameter values corresponding to these loops form the endpoints of the gaps in a Cantor set contained within the bifurcation diagram. Due to the presence of a Cantor set, there is a countable set of topologically non-equivalent germs of bifurcation diagrams even in generic one-parameter families, in contrast to bifurcations on the sphere.
\end{abstract}

\section{Introduction}

Suppose a vector field $v$ on a closed two-dimensional surface~$M$ has a hyperbolic saddle~$P$ with a separatrix loop~$\gamma$ which is orientation-preserving, that is, when we travel once along~$\gamma$, the local orientation does not change. Narrow neighborhoods of~$\gamma$ are cylinders then. The loop splits such a neighborhood into a monodromic and a non-monodromic semi-neighborhoods. In what follows, those separatrices of the saddle that are not involved in the loop will be called \emph{free}. If the loop~$\gamma$ is also non-contractible, it may happen that one of the free separatrices comes into a monodromic semi-neighborhood of the loop and winds onto the loop. If we want it to be the outgoing separatrix, we must require that the sum of eigenvalues of our saddle be non-positive. In fact, we will assume that it is strictly negative, since zero sum of eigenvalues \emph{and} the loop are not observed together in generic one-parameter families of vector fields. We will call saddles with negative sum of eigenvalues \emph{dissipative}.

\begin{figure}[h]
\center{\includegraphics[width=0.5\linewidth]{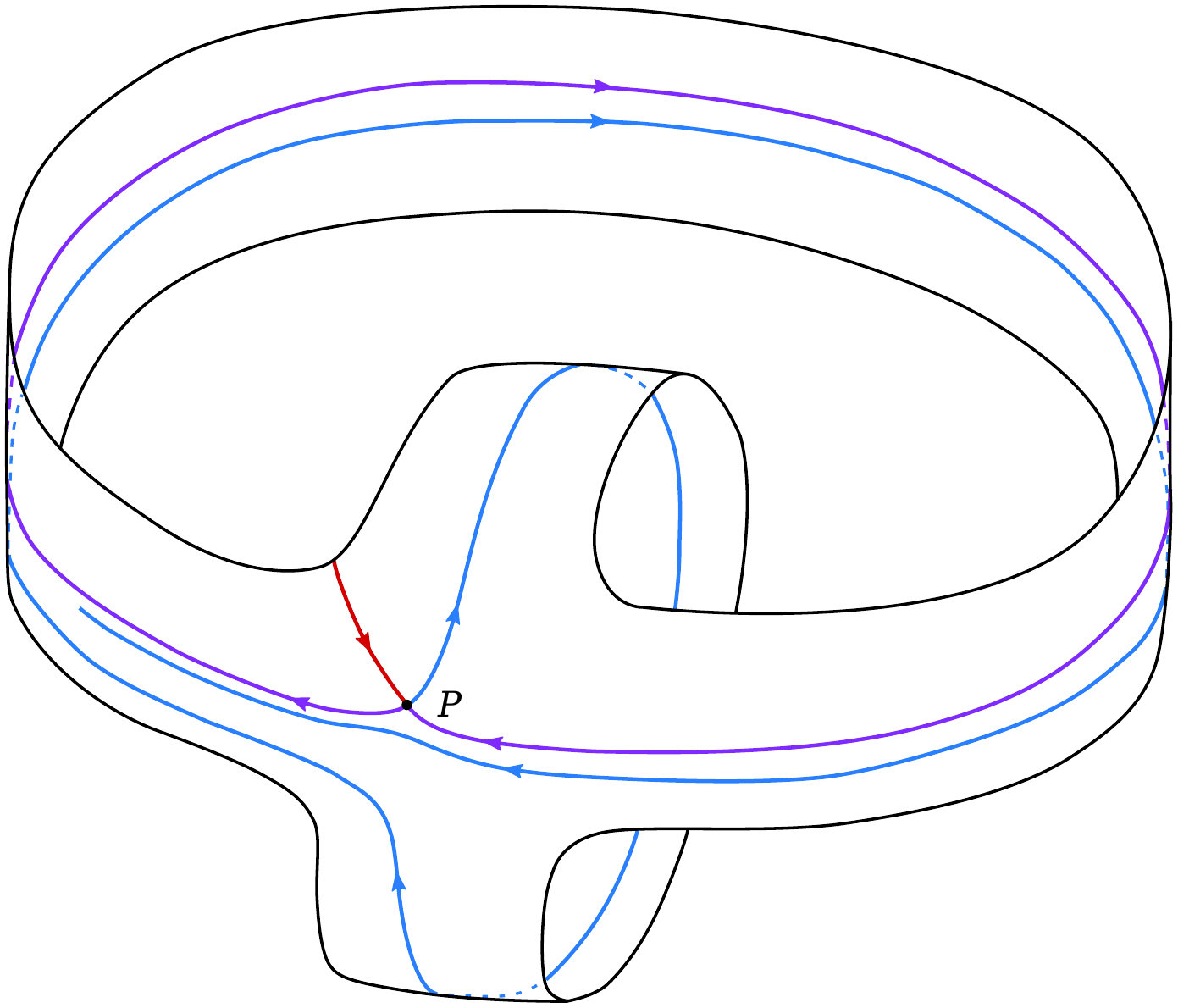}}
\caption{A narrow neighborhood of~$W^u(P)$ for~$v \in \orc(M)$.}\label{fig:0}
\end{figure}

There are two possibilities then: either the local orientation is preserved when we travel along the free separatrix back to the saddle~$P$, or it is inverted, and we will focus on the first case; the second will be addressed elsewhere. In the first case, a small neighborhood of the unstable manifold of the saddle is topologically a torus without a disk, also known as handle, see Fig.~\ref{fig:0}. There are no vector fields with such properties on the sphere, projective plane, and Klein bottle, because a handle cannot be embedded in these surfaces: when one cuts out an embedded handle and replaces it with a disk, this increases the Euler characteristic by 2 and preserves orientability, which leads to a contradiction for these three surfaces. On the other hand, it is not difficult to see (and we will show this below) that for all other closed connected surfaces such fields exist. We will need a shorthand notation for such fields.

\begin{defin} Let $v$ be a $C^\infty$-smooth vector field on a smooth closed two-dimensional surface~$M$. We say that~$v$ belongs to the class~$\orc(M)$ if the following two conditions hold.
\begin{enumerate}
\item The vector field $v$ has a dissipative saddle~$P$ with a saddle loop which is orientation preserving.
\item The unstable separatrix of~$P$ that is not involved in the loop winds onto it, and a narrow neighborhood of the unstable manifold of~$P$ is homeomorphic to a torus without a disk.
\end{enumerate}
\end{defin}

It is well known (see, e.g.,~\cite{IS} and references therein) what happens when a saddle loop of a vector field on the plane is unfolded in a generic one-parameter family. On one side of the critical parameter value a hyperbolic limit cycle is born. Moreover, if there are separatrices of other saddles that wind onto the loop, the critical value is accumulated from the other side by a sequence of parameter values that correspond to saddle connections between these saddles and the original one. These are called \emph{sparkling saddle connections}.

The case of a field $v \in \orc(M)$ is analogous, except that the connections are formed by the separatrices of the same saddle, so they can be called \emph{sparkling saddle loops}. It turns out that when a new loop of this type is formed, it is automatically accumulated by the free outgoing separatrix of the saddle, so unfolding this loop we get a new generation of loops and fields of class~$\orc(M)$, etc. It is this fascinating proliferation of sparkling loops that we want to draw the reader's attention to. We will see that the closure of the set of parameters that correspond to the presence of these loops is a Cantor set and these parameters themselves are the endpoints of the intervals of its complement. 

For a family of vector fields, the \emph{bifurcation diagram} is the set of parameter values that correspond to vector fields that are not structurally stable. Two bifurcation diagrams are called equivalent if the first diagram can be taken into the second one by a homeomorphism of the parameter space.
Any closed subset of the parameter space, regardless of its dimension, can be a bifurcation diagram for a sufficiently degenerate family: take a family obtained by multiplying some Morse-Smale field by a smooth function of the parameter that vanishes exactly at the closed subset of choice. However, it is natural to ask what diagrams, or rather different\footnote{Two germs of diagrams, say, at 0 are equivalent if they have equivalent representatives and the homeomorphism that realizes the equivalence of the representatives takes 0 to 0.} germs of diagrams, can occur in locally generic families.

It was conjectured by V. Arnold that for the case of vector fields on the sphere there exists but a finite number of pairwise non-equivalent germs of bifurcation diagrams that may occur in locally generic $k$-parameter families, for any~$k$ \cite{AAIS}.
This conjecture was disproved when a countable family of different germs of bifurcation diagrams was found in~\cite{KS} for three-parameter families that unfold the polycycle ensemble called \emph{the lips}. Then Yu.~Ilyashenko showed in \cite{I21} that infinitely many germs of bifurcation diagrams can be observed in two-parameter families that unfold a vector field that has a parabolic cycle accumulated by a separatrix of a saddle-node.
%A. Dukov noticed that such infinite family of germs had appeared, although not quite explicitly, in V.~Sh.~Roitenberg's PhD thesis~\cite{R}.
On the other hand, germs of generic one-parameter families of vector fields on the sphere have been classified in~\cite{IS},~\cite{IGS}, and~\cite{St}, and this classification yields, in particular, that these families admit only two nonempty germs of bifurcation diagrams up to equivalence: the first consists of one point at the origin, the second is the union of the point at the origin and a sequence that monotonically converges to it.
It turns out that due to the phenomenon discussed in the present paper, for one-parameter families on most other two-dimensional surfaces this is not the case.

%\begin{thm}\label{thm:diagv2}
%Let $M$ be a compact smooth two-dimensional surface other than the sphere, Klein bottle, and projective plane. Then locally generic one-parameter families of smooth vector fields on $M$ admit at least countably many non-equivalent germs of bifurcation diagrams.
%\end{thm}

\begin{mytheorem}\label{thm:diag}
Let $M$ be a smooth closed two-dimensional surface other than the sphere, Klein bottle, and projective plane. Then there exists a countable family of pairwise non-equivalent germs of bifurcation diagrams each of which is realized on an open set in the space of smooth one-parameter families of vector fields on~$M$.
\end{mytheorem}

%\begin{rem}\label{rem:struct}
%Let us denote the germs of bifurcation diagrams that will appear in the proof of Theorem~\ref{thm:diag} by~$K_n$, $n \in \nn$. 
%A representative of~$K_n$ has the following structure: it is a union of a Cantor subset of the real line, which is the closure of the set of parameter values that correspond to the presence of the sparkling loops, and a sequence of points. The Cantor subset contains zero and lies on one side from it. The sequence intersects each gap of the Cantor set (i.e., each connected component of the complement to the Cantor set in its convex hull) by exactly~$n$ points and has no other points.
%\end{rem}

We will prove Theorem~\ref{thm:diag} by modifying one-parameter families of vector fields on the torus studied by C.~Boyd in~\cite{B}. These families may be viewed as unfoldings of Cherry vector fields. Vector fields in these families have exactly two singularities, namely a saddle~$P$ and a sink~$\Omega$, admit a global Poincar{\'e} section, and their properties depend of the rotation number of the corresponding Poincar{\'e} map. In the parameter space, the bifurcation diagram of these families is a Cantor set whose gaps correspond to structurally stable vector fields. The endpoints of the gaps correspond to fields with a saddle loop; these fields will belong to the class $\orc(\TT)$ if we invert the time. In both cases, the rotation number of the Poincar{\'e} map is rational. The other points of the Cantor set correspond to the Cherry fields: the rotation number is irrational and the flow has a nontrivial recurrent trajectory.

We modify Boyd's families by adding $n$ additional Cherry cells to a neighborhood of the sink~$\Omega$. A Cherry cell is a pair of hyperbolic singularities: the first is a saddle and the second is a sink or a source that captures one of the separatrices of the saddle. Our Cherry cells will consist of a saddle and a source. For every Cherry cell and for every vector field of the family, both unstable separatrices of the saddle go to the sink~$\Omega$ and one stable separatrix goes to the source of the cell. The last stable separatrix leaves the neighborhood of the sink~$\Omega$ and can create a separatrix connection with the unstable separatrix of the original saddle~$P$, see Fig.~\ref{fig:V_k} below. These connections add~$n$ new bifurcation points into each gap of the Cantor set which was the original bifurcation diagram. The bifurcation diagrams for different integers~$n$ provide non-equivalent germs at points that correspond to the fields with a loop. The original Boyd families are known to be strongly structurally stable, which allows to prove that the new bifurcation diagram with additional points is observed in some open set of families. Finally, a small neighborhood of the sink can be cut out and replaced by an arbitrary two-dimensional surface with a circular boundary, which generalizes the result to the case of other surfaces.

In Theorem~\ref{thm:diag} different bifurcation diagrams appear due to the additional Cherry cells that sit far from the sparkling saddle loops. In Theorem~\ref{thm:main} we show that if we restrict our attention to a small neighborhood of the unstable manifold~$W^u(P)$, the bifurcations that happen when a field $v \in \orc(M)$ is perturbed will be essentially the same in any family.

%However, the main result of the present paper is the following theorem.

\begin{mytheorem}\label{thm:main}
Let $v$ be a vector field of class~$\orc(M)$. Then a generic $C^\infty$-smooth one-parameter family of vector fields~$V = \{v_\t\}_{\t \in [-1, 1]}$ unfolding~$v = v_0$ is strongly structurally stable in restriction to some small neighborhood of zero in the parameter space and some narrow neighborhood of~$W^u(P)$, the unstable manifold of the saddle~$P$ of~$v$, in the phase space. The bifurcation diagram of the restricted family is a Cantor set that contains zero and is located on one side of it; the endpoints of the gaps in the Cantor set correspond to the presence of ``sparkling'' separatrix loops for the saddle~$P$; whenever there is a separatrix loop for~$P$, the free unstable separatrix of~$P$ winds onto the loop. The bifurcation diagram has zero Hausdorff dimension. If~$W$ is another generic one-parameter family that unfolds a field $w\in\orc(N)$, then the families $V$ and $W$ are strongly topologically equivalent when restricted to some neighborhoods of the unstable manifolds of the corresponding saddles and some neighborhoods of zero in the parameter space.
\end{mytheorem}

%The last assertion of the theorem means that there is almost no hope that two Cantor sets of this origin can have robust intersection which could potentially lead to examples of generic one-parameter families where separatrix loops coexist with non-trivial recurrent trajectories for some values of the parameter.

To prove Theorem~\ref{thm:main}, we observe that $W^u(P)$ has a trapping neighborhood $U$, cut this neighborhood out, glue to it a disk with a source singularity, simultaneously for all parameter values close to 0, and invert the time. Some extra effort is needed to ensure that the resulting family of vector fields on the torus has Poincar{\'e} maps that are expansive (on the complement of a flat segment) and monotonic in the parameter. For such families the  strong topological equivalence was established by C.~Boyd, and the claim about the zero Hausdorff dimension follows from the results of J.~P.~P.~Veerman~\cite{V}.

Although we prove both theorems essentially by reducing them to the results obtained by C.~Boyd for families that unfold simple Cherry fields on~$\TT$, the reduction is not trivial. In particular, the monotonic dependence of the Poincar{\'e} maps on the parameter is easily lost when we perturb the family, and we prove that we can restore it by switching to a special chart on the transverse section.

It is worth noting that the complexity of the dynamics of Cherry fields, with their transitive quasi-minimal sets, originates, in a sense, from the unfolding of a saddle loop accumulated by the free separatrix in the orientation-preserving case.

%It seems, saddle loops accumulated by the free separatrix were hiding in plain sight in Boyd's families, shadowed by the irrational rotation numbers of the other fields nearby, whereas it would be natural to regard these sparkling loops, and not the rotation numbers, as the origin of bifurcations that happen in these families.

\section{Families of vector fields on surfaces}
In what follows, ``smooth'' always stands for ``$C^\infty$-smooth'', unless stated otherwise. The space of smooth vector fields on a smooth manifold~$M$ will be denoted by~$\mathrm{Vect}^\infty(M)$. %The $C^\infty$-smoothness assumption is used to refer to the results of Boyd, however, these results also hold for finitely smooth vector fields,
Our manifold usually will be a closed, connected two-dimensional surface. Such surface may be viewed, in a unique way, as a sphere, projective plane, or Klein bottle with zero or more handles attached. A \emph{handle} is a torus with a disk removed and by attaching a handle to a surface we mean taking the smooth connected sum of the surface and the handle. In what follows, a \emph{surface with a handle} is a closed, connected two-dimensional surface different from the sphere, projective plane, and Klein bottle. 

A \emph{family of vector fields} on a smooth manifold~$M$ with base $B$ is a smooth vector field on $B \times M$ that is tangent to the fibers $\{b\} \times M$. We will often refer to it as a map $V \colon B \to \mathrm{Vect}^\infty(M)$. From this point of view, a $C^1$-family of $C^\infty$ vector fields is a $C^1$-map from $B$ to $\mathrm{Vect}^\infty(M)$.
We will consider one-parameter families with base equal to a segment.

\begin{defin}
Two families $V = \{v_\t\}_{\t \in B}$ and $W = \{w_\tau\}_{\tau \in \hat{B}}$ of vector fields on homeomorphic manifolds~$M$ and~$\hat{M}$ are called \emph{strongly equivalent} if there exists a homeomorphism ${H\colon B\times M \to \hat{B}\times \hat{M}}$ of the form
$$(\t,x) \mapsto (h(\t), H_\t(x))$$
such that, for every $\t \in B$, the map $H_\t$ takes the phase portrait of $v_\t$ into the one of $w_{h(\t)}$, that is, takes trajectories to trajectories preserving the time-induced orientation. We say that the homeomorphism~$H$ \emph{realizes} this equivalence.
\end{defin}
A family $V$ is \emph{strongly structurally stable} if it is strongly equivalent to any family~$W$ that is sufficiently close in the space of families that we consider.
There are other notions of equivalence and hence stability for families; see, e.g.,~\cite[\S 1.1]{IKS} or~\cite{I21}.

%A \emph{local} family of vector fields is the germ at some point $b \in B$ of a global family. We will always assume that~$b = 0$ and so we will say that for a local one-parameter family the base is always~$(\RR, 0)$. Given a vector field $v$, a one-parameter unfolding of $v$ is a local family $V = \{v_\t\}_{\t \in (\RR, 0)}$ with $v_0 = v$.

%Two local families are called strongly equivalent if they have equivalent representatives and the first component (that is, $h$) of the map that realizes the equivalence of the representatives takes 0 to 0.

Throughout the paper we use standard dynamical systems notions and notation, for which the books~\cite{KH} and~\cite{PM} are a good reference.

\section{Proof of Theorem~\ref{thm:diag}}

%\subsection{Idea of the proof}
%C. Boyd~\cite{B} has described the bifurcation diagrams for a particular open set of one-parameter families of vector fields on the torus. The vector fields in these families lie in a vicinity of simple Cherry fields that have exactly two singularities, namely a saddle~$P$ and a sink~$\Omega$. The bifurcation diagrams of these families are Cantor sets. We modify the families considered by Boyd by adding to the basin of the sink, or rather, to a neighborhood of the sink, $n$ additional Cherry cells each of which consists of a saddle and a source, see Fig.~\ref{fig:V_k}. We can define a Cherry cell as a pair of hyperbolic singularities: the first is a saddle and the second is a sink or a source that captures one of the separatrices of the saddle. For every Cherry cell and for every vector field of the family, both unstable separatrices go to the sink~$\Omega$ and one stable separatrix goes to the source of the cell. The last stable separatrix leaves the neighborhood of the sink~$\Omega$ and can create a separatrix connection with the unstable separatrix of the original saddle~$P$. These connections add~$n$ new bifurcation points into each interval of the complement to the Cantor set which was the bifurcation diagram for the original Boyd family. The bifurcation diagrams for different integers~$n$ are not equivalent and, most importantly, they provide non-equivalent germs. Then it is not difficult to adapt this construction to the case of arbitrary surface with a handle.

\subsection{Boyd's families on the torus}
Consider the two-torus $\TT = \RR^2/\ZZ^2$. For a point in it, instead of writing $(x, y) + \ZZ^2$ or $[(x,y)]$, we will simply write~$(x, y)$. On the circle $\Sigma = \{(x,y)\in\TT \mid x = 0\}$ there is a natural coordinate~$y$ and the corresponding orientation. For a pair of points $a, b \in \Sigma$, we will denote by~$(a, b) \subset \Sigma$ the open arc that starts at $a$ and goes in the positive direction until it reaches~$b$. Whenever we differentiate a map from the circle to itself, it may be assumed that we take the derivative of the lift of our map to $\RR$, and if we differentiate a family of maps in the parameter, it may be assumed that the whole family is lifted to~$\RR$. At some point, we will need to change the coordinate on $\Sigma$. We will then assume that the new coordinate is given by some diffeomorphism $h: \RR/\ZZ \to \Sigma$, maybe parameter-dependent, and so there is a global ``circular chart''. We will also sometimes identify the points of $\Sigma$ and their vertical coordinates.

The following definition comes from~\cite{PM} and~\cite{B}.
\begin{defin}
We will say that a $C^\infty$-smooth vector field~$v$ on the torus~$\TT$ belongs to the class~$\N$ if the following conditions hold.

\begin{enumerate}
\item The vector field $v$ has exactly two singular points: a hyperbolic saddle~$P$ for which the sum of eigenvalues is positive and a hyperbolic sink~$\Omega$. 
\item The vector field~$v$ is transverse to the circle~$\Sigma = \{(x,y)\in\TT \mid x = 0\}$.
\item The local stable separatrices of the saddle~$P$ first cross~$\Sigma$ at points~$a$ and~$b$ (which depend on the field). One of the unstable separatrices of the saddle goes directly towards the sink~$\Omega$ without intersecting~$\Sigma$ and the second one first intersects~$\Sigma$ at a point~$c$.
\item For any point $y \in (a,b)\subset\Sigma$, the positive semi-orbit~$Orb_+(y)$ goes straight to the sink~$\Omega$ without re-intersecting~$\Sigma$ and for the points~$y\in(b,a) = \Sigma\setminus [a,b]$ the Poincar{\'e} map $f: (b, a) \to \Sigma$ is defined and is expansive: $f'(y) > \gamma > 1$ for all~$y$. 
\end{enumerate}
\end{defin}

\begin{rem}
\begin{itemize}
\item If a hyperbolic saddle of a vector field on a surface has eigenvalues of positive sum, such saddle is called anti-dissipative or area-expansive. For the vector field of class~$\N$ this property implies that $f'(y) \to +\infty$ as $y \to a-0$ or~$y \to b+0$.
\item The Poincar{\'e} map $f$ can be extended to the arc $[a,b]$ by setting $f([a,b]) = \{c\}$, where $c$ is the point where the unstable separatrix of~$P$ first intersects~$\Sigma$. Thus we get a continuous non-strictly monotonic map from~$\Sigma$ to itself, of degree one. This map is constant on~$[a, b]$, and in what follows we will refer to this segment as \emph{the flat segment} of~$f$. For the map~$f$ the rotation number $\rho(f) \in \RR/\ZZ$ is well-defined (see.~\cite[Sect. 11.1, p.~392]{KH}\footnote{See also \cite[Lemma 3 at p.~184]{PM} and \cite{Ha}, \cite{He}.}). In what follows, when we refer to the Poincar{\'e} map, we mean this extended map.
\item The class $\N$ is an open subset of the space~$\mathrm{Vect}^\infty(\TT)$ of smooth vector fields on the torus.
\end{itemize}
\end{rem}

The following theorem summarizes the results on the fields of class $\N$ obtained in~\cite{C}, \cite{PM}, and~\cite{B}; see also~\cite{KH}. In short, the theorem says that the dynamics can be described in terms of the rotation number.

\begin{thm}[\cite{C, PM, B}]\label{thm:BPM}
Let $v \in \N$ and let $f \colon \Si \to \Si$ be the corresponding extended Poincar{\'e} map.
\begin{enumerate}
\item If the rotation number~$\rho(f)$ is irrational, then
\begin{enumerate}
\item the free unstable separatrix of the saddle~$P$ intersects the transversal~$\Si$ infinitely many times, but it never intersects the closed arc~$[a, b] \subset \Si$;
\item the attraction basin $W^s(\Om)$ of the sink~$\Om$ is dense in~$\TT$;
\item its complement $\TT \setminus W^s(\Om)$ is a transitive quasi-minimal set\footnote{A quasi-minimal set is the closure of a trajectory that is recurrent in both directions. It is often also assumed that this recurrent trajectory is nontrivial (i.e., neither a singular point nor a cycle)~\cite[p.~83]{ABZh}.} of zero measure that intersects~$\Sigma$ by a Cantor set of zero measure\footnote{Moreover, this Cantor set has zero Hausdorff dimension~\cite{V}.}.
\end{enumerate}
If for the fields $v, w \in \N$ the Poincar{\'e} maps have the same irrational rotation numbers, these fields are orbitally topologically equivalent.

\item If the rotation number $\rho(f)$ of the Poincar{\'e} map for the vector field~$v\in\N$ is rational, then there are two possible cases:
  \begin{enumerate}
  \item either the free unstable separatrix of the saddle~$P$ intersects the transversal~$\Si$ finitely many times, the last intersection being at the point~$a$ or~$b$, and thus, there is a separatrix loop;
  \item or the free unstable separatrix of the saddle~$P$ intersects~$\Si$ finitely many times and the last intersection is inside the open arc~$(a, b) \subset \Si$; then~$v$ is Morse-Smale.
  \end{enumerate}
Two fields of type (a) (resp. (b)) with the same rotation number are orbitally topologically equivalent.
\end{enumerate}
\end{thm}

Note that when the rotation number is irrational, the field cannot be Morse-Smale due to the presence of a (non-trivial) quasi-minimal set, which is part of the non-wandering set (alternatively, we may argue that the rotation number can be made rational by a small perturbation).

Colin Boyd proved a strong stability result for quite specific families of vector fields of class~$\N$.

\begin{thm}[C. Boyd, \cite{B}]\label{thm:Boyd}
Let $V = \{v_\t\}_{\t \in [0, 1]},\; v_\t \in \N$, be a $C^1$-smooth one-parameter family of vector fields such that
$$f_\t(\cdot) = f_0(\cdot) + \t,$$
where $f_\t$ is the Poincar{\'e} map for the field~$v_\t$. If~$v_0$ is Morse-Smale, then the family~$V$ is strongly structurally stable in the space of $C^1$-paths in~$\N$.
\end{thm}

Note that the family $V$ in this theorem is very degenerate; in particular for all Poincar{\'e} maps the flat segment is the same and, as it is easy to see when looking at the asymptotic of the Poincar{\'e} map at the points~$a, b$, the ratio of eigenvalues of the saddle in this family does not depend on the parameter. Nevertheless, by the theorem, there are $C^1$-families of general position which are equivalent to it, and the same can be said about $C^\infty$-families.

%%%%%%%%%%%%%%%%%%%%%%%%%%%%%%%%%%%%%%%%%%%%%%%%%%

\subsection{Bifurcation diagrams of Boyd's families}\label{sect:bifdiagboyd}

Boyd remarks~\cite[p. 28]{B} that it is easy to see that the bifurcation diagrams for families from Theorem~\ref{thm:Boyd} are Cantor sets. He does not give the proof, so we present this argument for completeness in this section.

Consider a family~$V$ as in Theorem~\ref{thm:Boyd}. The existence of such a family is obvious: one can take a vector field~$v \in \N$ and start rotating it in a vertical strip near~$\Sigma$, on the ``incoming'' side; an explicit construction can be found in~\cite[p.~186]{PM} (see also~\cite[p.~464]{KH}).

\begin{prop}\label{prop:bdiag}
Let $V = \{v_\t\}_{\t \in [0, 1]},\; v_\t \in \N$, be a $C^1$-smooth one-parameter family of vector fields such that
$$f_\t(\cdot) = f_0(\cdot) + \t$$
and $v_0$ is Morse-Smale. Then the bifurcation diagram of $V$ is a Cantor set that coincides with the closure of the set of parameters that correspond to the presence of separatrix loops. The same holds for any other family in a small neighborhood of~$V$ in the space of $C^1$-families of $C^\infty$ vector fields.
\end{prop}

\begin{proof}

Consider a neighborhood~$\F$ of the family~$V$ in the space of $C^1$-families of $C^\infty$-smooth vector fields. Let the neighborhood $\F$ be so small that all families in~$\F$ are strongly equivalent to each other and contain only vector fields of class~$\N$. Then the bifurcation diagrams of these families are also equivalent to each other: for every such family a parameter value is not in the bifurcation diagram if and only if the corresponding vector field is Morse-Smale, and the homeomorphism that realizes the strong equivalence of two families takes Morse-Smale fields to Morse-Smale ones.\footnote{With our definition of bifurcation diagrams, it is easy to construct examples of strongly equivalent families with non-equivalent bifurcation diagrams. E.g., imagine a family with a hyperbolic cycle and an equivalent family with a corresponding cycle that is not hyperbolic for some isolated parameter value: this value will be in the diagram for the second family. So, bifurcation diagram is not an invariant of topological classification of families. It becomes an invariant if we consider only families that do not contain vector fields that are not structurally stable, but are orbitally topologically equivalent to structurally stable ones. Had we defined bifurcation diagram as the set of parameter values that do not have a neighborhood where all corresponding fields are equivalent, it would automatically be an invariant. However, we prefer our bifurcation diagrams to reflect the lack of structural stability rather than to indicate that the bifurcation is truly observed in the family.}

Let us look at the bifurcation diagram~$\Bif(V)$ of the family~$V$. Consider the corresponding family~$F = \{f_\t\}$ of Poincar{\'e} maps and the function~$r(\t) = \rho(f_\t)$. By Theorem~\ref{thm:BPM}, if $r(\t)$ is irrational, the corresponding vector field is not Morse-Smale, which yields that $A = \{\t\colon r(\t)\notin\QQ/\ZZ\} \subset \Bif(V)$. Furthermore, the family~$F$ is monotonic in~$\t$; therefore~$r$ is non-decreasing and, moreover, it is strictly increasing at points where it has irrational values, see~\cite[Prop. 11.1.8-9]{KH}, so a fixed irrational value is attained at isolated points. Since $\t \in [0, 1]$ and condition $f_\t = f_0 + \t$ holds, any irrational value is attained at exactly one point.

Recall that we denote by $c$ the point of the first intersection of the unstable separatrix of $P_\t$ with $\Sigma$; this point depends on $\t$ now: $c = c(\t)$. Fix some~$\t$ such that~$r(\t)$ is rational. By Theorem~\ref{thm:BPM}, the unstable separatrix of~$P_\t$ intersects the arc~$[a, b]\subset\Sigma$. This means that for some~$k\in\nn$ we have
$f_\t^{k-1}(c) \in [a,b].$
Since $f_\t(x)$ is non-decreasing in $x$ and strictly increasing in~$\t$, we can conclude that $f_\t^{k-1}(c)$ is also strictly increasing in~$\t$. This means that for the chosen rotation number the separatrix crosses $[a, b]$ when~$\t$ belongs to some segment. The endpoints of the segment correspond to the presence of loops (the unstable separatrix comes to~$a$ or~$b$), whereas the interior of the segment corresponds to Morse-Smale vector fields. In the neighborhood of the endpoints the rotation number is non-constant. Indeed, when $f_\t^{k-1}(c)$ goes past, say, $b$ as~$\t$ increases, the point~$c$ can no longer be periodic with period~$k$, and there are no other periodic orbits of that period, because their existence would imply, by Theorem~\ref{thm:BPM}, that the field either is Morse-Smale or has a loop, but in both cases the point~$c$ would be periodic with the same period~$k$. Thus, the rotation number has to change when we go past the endpoint.

Since any $\t$ with irrational $r(\t)$ is approximated by segments where $r(\cdot)$ is rational and by their endpoints that correspond to loops, we conclude that
$$\Bif(V) = \overline{\{\t\colon r(\t)\notin\QQ/\ZZ\}} = \overline{\{\t\colon v_\t  \text{ has a loop}\}}.$$
Now it is clear that $\Bif(V)$ is perfect and nowhere dense and hence it is a Cantor set. The same can be said about any family in the neighborhood~$\F$.
\end{proof}

%%%%%%%%%%%%%%%%%%%%%%%%%%%%%%%%%%%%%%%%%%%%%%%%%%

\subsection{Adding Cherry cells}
Consider a family~$V$ as in Proposition~\ref{prop:bdiag}, but also~$C^\infty$-smooth. Moreover, assume for convenience that outside of some thin vertical strip that contains no singularities and is on the ``incoming'' side of the circle~$\Sigma$ the vector fields of our family are exactly the same for all parameter values.
 
We want to perform, on the family~$V$, a surgery that will add a Cherry cell into a neighborhood of the sink~$\Omega$. 
Near the sink the vector fields of our family~$V$ do not depend on the parameter. Consider a small flow box~$\Pi_1$ near the sink. Choose the flow-box so that it have no intersection with the unstable separatrix $u$ of~$P$ that always goes to~$\Omega$. Switch to the rectifying coordinates of this flow box, and for each parameter value replace the constant flow inside the box with a field that has a Cherry cell with a saddle~$P_1$ and a source~$A_1$. The outgoing separatrices of~$P_1$ have to go to the sink~$\Omega$, one incoming separatrix is captured by the source~$A_1$ and the other is, in a sense, free, and has to cross the transversal~$\Sigma$, see Fig.~\ref{fig:V_k}.

As we do this surgery, we can also make sure that, first, near the boundary of the flow box the fields of the family remain the same, and second, the Cherry cell slowly moves ``downwards'' as the parameter increases: that is, we want the point~$d$ of intersection between the free stable separatrix of~$P_1$ and~$\Sigma$ to be strongly monotonic in~$\t$ in the sense that we must have $d'(\t) < -\beta <0$ for all~$\t$.\footnote{Here we denote the coordinate of the point~$d$ on the vertical circle by the same symbol.} To achieve that we glue in the ``same'' Cherry cell for different~$\t$, but we shift its position by~$\delta \cdot \t$ in the direction perpendicular to the constant field that we originally have in the rectifying coordinates on~$\Pi_1$; here~$\delta$ is a small constant.

Denote by $B_1$ a small disk that is independent of the parameter, contains a neighborhood of the flow box~$\Pi_1$, does not intersect with the separatrix~$u$, and for every parameter value is contained in the basin of the sink~$\Om$. Denote by~$V_1$ the special family obtained from~$V$ by adding one Cherry cell as described, and denote a small neighborhood of~$V_1$ in the space of $C^\infty$-families by~$\F_1$.
We take $\F_1$ so small that for any family $W\in~\F_1$ we have the following:
\begin{enumerate}
\item the continuation of the saddle $P_1$ and source $A_1$ is in~$B_1$ for every $\t\in[0, 1]$;
\item there exists a smoothly parameter-dependent curvilinear rectangle $\tilde{\Pi}_1(\t)$ such that
\begin{itemize}
\item $P_1(\t), A_1(\t) \in \tilde{\Pi}_1(\t) \subset B_1$;
\item two opposite edges of $\tilde{\Pi}_1(\t)$ are transverse segments and the other two edges are segments of trajectories;
\item inside $\tilde{\Pi}_1(\t)$ the field can be replaced by a field that is smoothly equivalent to a constant one in such a way that this yields a family of class~$\F$;
\end{itemize}
\item\label{pr:monot} for the family $W$, the point of intersection between the free stable separatrix of the saddle~$P_1$ and the circle~$\Sigma$ strongly monotonically depends on the parameter in the same sense as above: namely, the derivative of its vertical coordinate in $\t$ is negative.\footnote{Here and below it is important that for two families of vector fields which are close to each other the families of local stable (or unstable) manifolds of the hyperbolic continuations of some saddle are also close. Here we only need them to be $C^1$-close, but below we will need $C^3$-closeness. However, this also holds in~$C^r$, $1 \le r < \infty$. This can be justified in the following way. Consider the family of vector fields with a parameter-dependent saddle $P_\t$ as one vector field $v_\t\frac{\partial}{\partial \overrightarrow{\boldsymbol{x}}} + 0\cdot \frac{\partial}{\partial \t}$ on~$\TT\times [0, 1]$. For this field the local central-stable manifolds of the points~$(P_\t, \t)$ are uniquely defined and coincide with each other. Most importantly, these manifolds are $C^r$-smooth and continuously depend on the field in the $C^r$-topology, see~\cite[Theorems 5.1, 5A.1]{HPS}. In particular, if we consider the point of intersection of a stable separatrix of $P_\t$ with a fixed smooth transversal as a function of $\t$, a $C^r$-perturbation of the family will yield a~$C^r$-close function.} %The two families of stable manifolds of the points~$P_\t$ are obtained as transverse intersections between these central-stable manifolds, which are close, and planes of the form $\{\t = const\}$; therefore these families are~$C^r$-close as well.}
\end{enumerate}

\begin{figure}
  \begin{center}

  \begin{tikzpicture}[scale=1.6]
  
  % grid
  %\draw[step=0.2cm,gray,very thin] (0,0) grid (4,4);

  %coordinates
  \coordinate (P) at (2.4, 2.2);
  \coordinate (u) at (2.05, 2.1);
  \coordinate (a) at (0, 0.7);
  \coordinate (b) at (0, 2.8);
  \coordinate (c) at (4, 2.5);
  \coordinate (Om) at (1.65, 2.0);
  \coordinate (P1) at (0.65, 2.2);
  \coordinate (P2) at (0.7, 1.35);
  \coordinate (A1) at (1.03, 2.22);
  \coordinate (A2) at (1.07, 1.44);

  %separatrices
  % for P
  \draw[blue, very thick] plot [smooth, tension = 0.7] coordinates{(Om) (P) (c)};
  \draw[red, very thick]  (b) .. controls (1.3, 2.8) and (2.3, 2.7) .. (P) .. controls (2.5, 1.9) and (1.9, 1.) .. (a);
  % for P_1 and P_2
  \draw[red, very thick] plot [smooth, tension = 0.7] coordinates{(0, 2.12) (0.33, 2.17) (P1) (A1)};
  %\draw[blue, very thick]  (1.29, 2.32) .. controls (1.1, 2.45) and (0.7, 2.35) .. (P1) .. controls (0.7, 2.03) and (1.1, 2.1) .. (1.21, 2.0);
  % W^u(P1)
  \draw[blue, very thick]  (Om) .. controls (1.44, 1.92) and (1.3, 2.27) .. (1.21, 2.30) .. controls (1.1, 2.4) and (0.68, 2.4) .. (P1) .. controls (0.68, 2.0) and (1.05, 2.085) .. (1.1, 2.065) .. controls (1.28, 2.085) and (1.42, 1.9) .. (Om);

  \draw[red, very thick] plot [smooth, tension = 0.7] coordinates{(0, 1.25) (0.31, 1.29) (P2) (A2)};
  %\draw[blue, very thick]  (1.29, 1.7) .. controls (1.1, 1.53) and (0.72, 1.57) .. (P2) .. controls (0.7, 1.2) and (1.15, 1.22) .. (1.47, 1.58);
  \draw[blue, very thick]  (Om) .. controls (1.44, 1.94) and (1.35, 1.72) .. (1.29, 1.7) .. controls (1.1, 1.53) and (0.72, 1.57) .. (P2) .. controls (0.7, 1.2) and (1.15, 1.22) .. (1.43, 1.59) .. controls (1.56, 1.79) and (1.44, 1.9) .. (Om);

  %arrows on separatrices
  % for P
  \draw[->, red, thick] (2.0, 2.59) -- (2.07, 2.55);
  \draw[->, red, thick] (2.05, 1.54) -- (2.1, 1.59);
  \draw[->, blue, thick] (2.85, 2.29) -- (2.95, 2.31);
  \draw[<-, blue, thick] (2.0, 2.1) -- (2.05, 2.115);

  % for P_1 and P_2
  \draw[->, red, thick] (0.5, 2.187) -- (0.55, 2.194);
  \draw[->, red, thick] (0.5, 1.315) -- (0.55, 1.322);
  \draw[->, red, thick] (0.81, 2.207) -- (0.80, 2.206);
  \draw[->, red, thick] (0.83, 1.382) -- (0.82, 1.38);

  %\draw[<-, blue, thick] (0.85, 2.345) -- (0.84, 2.342);
  \draw[<-, blue, thick] (0.98, 2.367) -- (0.97, 2.366);
  %\draw[<-, blue, thick] (0.85, 2.082) -- (0.84, 2.084);
  \draw[<-, blue, thick] (0.93, 2.069) -- (0.91, 2.069);
  \draw[->, blue, thick] (0.86, 1.515) -- (0.88, 1.524);
  \draw[->, blue, thick] (0.88, 1.261) -- (0.9, 1.260);

  %saddle P and points a, b, c
  \filldraw [color = black] (P) circle (1pt) node[above right] {$P$};
  \filldraw [color = black] (a) circle (0.5pt) node[left]  {$a$};
  \filldraw [color = black] (b) circle (0.5pt) node[left]  {$b$};
  \filldraw [color = black] (c) circle (0.5pt) node[right] {$c$};
  \filldraw [color = red] (Om) circle (1pt) node[above, color = black] {$\Omega$};
  \filldraw [color = black] (P1) circle (0.5pt);
  \node[above left] at (0.73, 2.2)  {{\scriptsize $P_1$}};
  \filldraw [color = black] (P2) circle (0.5pt);
  \node[below left] at (0.75, 1.35) {{\scriptsize $P_2$}};
  \filldraw [color = blue] (A1) circle (0.5pt);
  \filldraw [color = blue] (A2) circle (0.5pt);

  %torus label
  \node[above] at (0, 4) {$\mathbb{T}^2$};

  %separatrix label
  \node[below] at (u) {$u$};

  %Transversal label
  \node[below right] at (0, 4) {$\Sigma$};
  \node[below left]  at (4, 4) {$\Sigma$};

  %torus rectangle
  \draw[thick] (0,0) rectangle (4,4);
  \end{tikzpicture}
  \caption{Fields of the family~$V_k$ for $k = 2$.}\label{fig:V_k}
  \end{center}
\end{figure}

We define special families $V_k$ and their open neighborhoods~$\F_k$, $k\in\nn$, analogously, but in $\F_k$ families have~$k$ Cherry cells instead of one and these cells are ``independent'' in the sense that for all parameter values the unstable separatrices of the saddles go to the sink~$\Omega$ and for each saddle the free stable separatrix intersects the transversal~$\Si$. The conditions above must hold for each Cherry cell, with disjoint $B_j\supset \tilde{\Pi}_j(\t)$.

\subsection{Required bifurcation diagrams for families of fields on~\texorpdfstring{$\TT$}{T2}}

In this section we will show that for the families from open sets~$\F_k$ the bifurcation diagrams have the required structure, that is, these are Cantor sets with~$k$ isolated points added into each gap.

Note that, by construction of the set~$\F_k$, for any family~$W \in \F_k$ there exist (parameter-independent) neighborhoods $B_1, \dots, B_k$ of the Cherry cells and a family $\overline{W} \in \F$ such that~$W$ and~$\overline{W}$ coincide in restriction to~$\TT\setminus\sqcup B_j$.
This implies that~$\Bif(\overline{W}) \subset \Bif(W)$. Indeed, if the field $\overline{w}_\t$ has a saddle loop, the same holds for the field~$w_\t$, and $\Bif(\overline{W})$ is the closure of the set of such parameter values. Now, consider a parameter value~$\t\notin \Bif(\overline{W})$. For this value the vector field~$\overline{w}_\t$ is Morse-Smale. The field~$w_\t$ cannot have cycles that intersect the disks~$B_j$, so all its cycles are the same as for the field~$w_\t$. Actually, the Poincar{\'e} map being expansive implies that the field has only one repelling cycle, and it is hyperbolic. All singular points of~$w_\t$ are hyperbolic as well. The non-contractible repelling cycle cuts the torus into a cylinder, which admits no nontrivial recurrence. Hence, the non-wandering set of~$w_\t$ contains only cycles and singularities, since there are no saddle loops for the saddle~$P$. Peixoto's theorem on structural stability yields then that the only way the field~$w_\t$ can be not structurally stable is by having a saddle connection between the saddle~$P$ and some saddle~$P_j$ (connections between saddles~$P_i$ and~$P_j$, including loops, are impossible by construction).

We want to show that in each interval of the complement to~$\Bif(\overline{W})$ there is exactly one parameter value that corresponds to the field with a separatrix connection between~$P$ and~$P_j$, for each~$j$.
For that we need the following proposition which we will prove below in section~\ref{sect:monoproof}.

\begin{prop}\label{prop:monot}
Let $V= \{v_\t\}_{\t \in [0, 1]}, \; v_\t \in \N$, be the $C^\infty$-family defined above in Section~\ref{sect:bifdiagboyd}, and let some $\e > 0$ be fixed. Then, if a family~$W$ is sufficiently close in the $C^\infty$-topology to the family~$V$, there exists a ($C^3$-smooth, at least) parameter-dependent coordinate on~$\Si$ that coincides with the original coordinate outside $\e$-neighborhoods of the points~$a$, $b$ and such that the family of the Poincar{\'e} maps $g_\t$ of the family~$W$, when written in the new coordinate, is strongly monotonic in the parameter~$\t$:
$$\forall y\in\RR/\ZZ,\; \forall \t_0 \in [0, 1], \text{ we have } \left.\ddt\right|_{\t = \t_0}\hat{g}_{\t}(y) > 0,$$
where $\hat{g}_\t$ is $g_\t$ written in the new coordinate.
\end{prop}

\begin{rem}\label{rem:fail}
Without switching to parameter-dependent coordinates, this would not hold. Consider, for example, a family obtained from~$V$ by precomposition with a vertical shift by $\e \t$:
$$w_\t(x,y) = v_\t(x, y - \e \t).$$
If the constant $\e > 0$ is small, this family is close to the family~$V$, but the derivative of the Poincar{\'e} map $g_\t(y) = f_\t(y - \e \t) + \e \t$ in $\t$ tends to $-\infty$ as $y\to a + \e\t - 0$ or $y\to b + \e\t + 0$, as direct calculation shows. The reason is that we multiply by~$(f_\t)'_y$, which is unbounded. However, it is clear that if we look at this family of Poincar{\'e} maps via the parameter-dependent coordinates change $(x,y)\mapsto (x, y - \e \t)$, we will again see the family~$\{f_\t\}$ that is monotonic in the parameter.
\end{rem}

Let us continue to prove the theorem. Assume that for some fixed parameter value there is a saddle connection between the saddle~$P$ and, say, the saddle~$P_1$.
Proposition~\ref{prop:monot} implies that in an appropriate chart the iterates~$\hat{g}_\t^{\circ n}(\cdot)$ of the Poincar{\'e} map are monotonic in the parameter. Take the last (if we count from the saddle~$P$) point of intersection between the unstable separatrix of $P$ that forms the connection and the transversal~$\Sigma$. This point monotonically depends on the parameter~$\t$. On the other hand, by property~\ref{pr:monot} of the class~$\F_k$, the vertical coordinate  of the first intersection point between the free stable separatrix of~$P_1$ with this circle monotonically decreases with the parameter. Note that we can assume that this point is always far from the points~$a, b$, so monotonicity is preserved when switching to the new chart. Therefore, the saddle connection happens at a unique point in the interval of the complement to $\Bif(\overline{W})$. The same argument works for connections between~$P$ and other saddles, hence the bifurcation diagram~$\Bif(W)$ consists of the Cantor set $K = \Bif(\overline{W})$ and a countable set of points that intersects every interval of the complement to~$K$ by exactly~$k$ points which correspond to separatrix connections between the saddle~$P$ and the saddles~$P_1, \dots, P_k$.

Since we proved this for an arbitrary $W \in \F_k$, we now have a countable family of open sets of families with different bifurcation diagrams that have the required structure. Since the germs at zero of these diagrams are also different, we proved Theorem~\ref{thm:diag} for the case of the torus (modulo Proposition~\ref{prop:monot}). It is not difficult to deduce from Boyd's Theorem~\ref{thm:Boyd} that in each $\F_k$ the families are strongly equivalent to each other, but we will not need this. 
In the following section we adapt this construction to an arbitrary surface with a handle.

\subsection{The case of arbitrary surface \texorpdfstring{$M$}{M}}
In this section we generalize the construction above by replacing the neighborhood of the sink~$\Omega$ by a surface with a Morse-Smale vector field on it.

An arbitrary smooth closed surface $M$ with a handle can be obtained from another surface~$N$ by gluing a handle: $M = N \# \TT$. For the surface~$N$, there exists a Morse-Smale vector field with a hyperbolic sink\footnote{Indeed, any vector field on~$N$ with a hyperbolic sink can be $C^1$-perturbed into a Morse-Smale field, and the latter can be approximated by a $C^\infty$-smooth Morse-Smale field which still has the sink.} and, therefore, there is also a Morse-Smale vector field $v_{N}$ with a small contractible hyperbolic attracting cycle that bounds a disk with a single singularity --- a hyperbolic source.
%The second field can be obtained from the first one by changing it in a vicinity of the sink in such a way that it undergoes the Andronov-Hopf bifurcation, or, alternatively, we can cut out a disk that contains the sink and glue in a disk with a vector field that has the required cycle and the source.
Since $v_{N}$ is Morse-Smale, it has a neighborhood~$\U\ni v_{N}$ where all vector fields are orbitally topologically equivalent.

Fix some $k$ and consider the special family $V_k \in \F_k$. Recall that its sink~$\Omega$ does not depend on the parameter, draw a small transverse circle around~$\Omega$, cut out the disk~$D\ni\Omega$ bounded by the circle (the disk must be small and should not intersect the disks~$B_j$), and denote what is left of the torus, namely $\TT\setminus D$, by~$T_0$. Note that by construction the restrictions of the fields of our family to the vicinity of $\partial T_0$ do not depend on the parameter.

Consider a constant family $\hat{V}$ whose vector fields coincide with~$v_{N}$ for all values of the parameter. For the vector field~$v_{N}$, draw a small transverse circle around the aforementioned source which is inside the attracting cycle, cut out the disk bounded by the circle and denote the rest of the surface by~$N_0$. Now let us smoothly glue $T_0$ to $N_0$ along the neighborhoods of the boundaries in such a way that on the resulting surface (diffeomorphic to~$M$) we get a smooth family of vector fields~$V_{k,M}$ --- the result of ``gluing'' the families~$V_k$ and~$\hat{V}$ together.

Note that for any family $V_{M}$ that is sufficiently close to~$V_{k, M}$ there exists a family~$V_{\TT} \in \F_k$ such that the restrictions of~$V_{M}$ and~$V_{\TT}$ on~$T_0$ coincide.\footnote{Here we assume that $T_0$ is simultaneously a subset of the torus and the surface~$M$ and $N_0$ is a subset of~$M$ and~$N$.}
Moreover, if the family~$V_{M}$ is sufficiently close to~$V_{k, M}$, the restriction $V_{M}|_{N_0}$ is close to the restriction~$\hat{V}|_{N_0}$, i.e., it is almost constant, and so every vector field of the restricted family coincides with the restriction to~$N_0$ of some field from the neighborhood~$\U\ni v_{N}$ where all vector fields are topologically equivalent.

The families $V_M$ and~$V_{\TT}$ have the same bifurcation diagrams. Indeed, the inclusion $\Bif(V_\TT) \subset \Bif(V_M)$ is obvious. On the other hand, if for some $\t_0$ the field $v_{\TT,\t_0}$ is Morse-Smale, we immediately have that the field~$v_{M,\t_0}$ is Morse-Smale as well, which yields the opposite inclusion.

Since this argument works for arbitrary $k\in\nn$, we get a countable family of open sets~$\F_{k, M} \ni V_{k, M}$ of one-parameter families of vector fields that have the required bifurcation diagrams.

\subsection{Germs of bifurcation diagrams}
For a family from the set~$\F_{k, M}$, the germs of bifurcation diagrams at parameter values that correspond to vector fields with separatrix loops for the saddle~$P$ are all equivalent. These germs have the form of a germ of a Cantor set on one side of the critical parameter value with~$k$ points added into each gap. For two families taken from the sets~$\F_{k, M}$ and~$\F_{j, M}$, $k \ne j$, these germs are not equivalent; they are distinguished by the number of isolated points in the gaps of the Cantor set. This proves Theorem~\ref{thm:diag} modulo Proposition~\ref{prop:monot}.
\begin{rem}
The same applies to germs at parameter values where the rotation number of the extended Poincar{\'e} map is irrational. Recall, however, that the equivalence of germs is established by a germ of homeomorphism that takes the selected point to the selected point, so this formally gives us another countable family of non-equivalent germs. For a fixed $k$, the germs at different parameter values of this type are the same, since any point of a Cantor set that is not an endpoint of a gap can be taken to any other such point by a homeomorphism of $\RR$ that takes this Cantor set into itself.
\end{rem}

\subsection{Proof of Proposition~\ref{prop:monot}}\label{sect:monoproof}
We wish we could simply say that, if two families of vector fields are close, they have families of Poincaré maps which are close, and, since for the family~$V$ the derivative of the Poincar{\'e} map in the parameter is always positive, the same must hold for a family~$W$ close to~$V$. However, we deal with Poincar{\'e} maps that were extended through the singularity, and therefore we cannot reason like that. Also recall Remark~\ref{rem:fail}.

For the special family~$V$, the points where the stable separatrices of the saddle~$P_\t$ intersect the circle~$\Sigma$ do not depend on~$\t$. Therefore, for a family~$W$ close to~$V$ these points change only slightly and with almost zero speed when the parameter varies, so there exists a parameter-dependent change of coordinates on the torus with the following properties. The coordinate change must be at least $C^3$-smooth, with uniformly small derivative in the parameter, uniformly $C^1$-close to the identity for all $\theta$, equal to the identity outside small, but chosen beforehand neighborhoods of the points $a, b$, and, moreover, in the new chart these points of intersection with the invariant manifolds must appear to be independent of the parameter.\footnote{E.g., we can take the coordinates change
$$y \mapsto \varphi_1(y)y + \varphi_2(y)(y + a - a_\t) + \varphi_3(y)(y + b - b_\t)$$
where $a_\t, b_\t$ are the original parameter-dependent coordinates of the intersection points between the local stable separatrices of~$P$ and the transverse circle~$\Sigma$ for the family~$W$ and $\{\varphi_j\}$ is a partition of unity on~$\Si$ such that~$\varphi_2$ and~$\varphi_3$ have supports in $\e$-neighborhoods of~$a, b,$ respectively, and are equal to 1 in the corresponding $\e/2$-neighborhoods. Note that $a-a_\theta, b-b_\theta, (a_\theta)'_\theta, (b_\theta)'_\theta$ are small when $W$ is close to~$V$.} Since $\ddt f_\t = 1$,  we can assume that $\ddt \hat{f}_\t > 1/2$. Here and below the ``hat'' indicates that we are looking at something in our new coordinates. 

The new coordinates provide a ``circular chart'' on $\Sigma$, that is, a diffeomorphism between~$\Sigma$ and $S^1 = \RR/\ZZ$. We will denote by~$a, b$ the points in $S^1$ that correspond to the points of intersection with the separatrices. Recall that~$\{g_\t\}$ are the Poincar\'e maps for~$W$. All points $y \in [a_\t,b_\t] \subset \Sigma$ have the same image under~$g_\t$ that coincides with the point $c_\t$ of intersection between the local unstable separatrix of~$P_\t$ with~$\Sigma$. This point depends at least $C^3$-smoothly on the parameter, therefore the derivative $\left.\frac{\partial}{\partial \t}\right|_{\t = \t_0} \hat{g}_\t$ is defined for the corresponding points in $[a, b] \subset S^1$, coincides with $\hat{c}'(\t_0)$, and is positive for all~$\t_0$, provided that~$W$ is sufficiently close to~$V$ and the coordinate change is sufficiently close to the identity and has small derivative in~$\t$. For a point $y \in S^1\setminus [a, b]$ this derivative is also defined and positive for all~$\t_0$ if~$W$ is close to~$V$. %However, if $\frac{\partial}{\partial \t}|_{\t = \t_0} \hat{g}_\t(y)$ were not continuous at the points~$a$ and~$b$, it could in principle turn out that for different points $y \in S^1 \setminus [a, b]$ this derivative is positive for the families in different, decreasing neighborhoods of the family~$V$ and it could be possible to approximate~$V$ by families that have negative derivative in the parameter at some points. [Remark~\ref{rem:fail} makes this redundant.]

We will show that any $\t_0$ has a neighborhood where, as $y \to b + 0$, we have the uniform convergence
\be\label{eq:conv}
\ddt {\hat{g}}_\t(y) \rightrightarrows \ddt{\hat{g}}_\t(b) = \hat{c}'(\t) > 0.
\ee
This will allow  us to take a finite cover of the compact parameter space by these neighborhoods and then choose $\alpha > 0$ such that $\left.\ddt\right|_{\t = \t_0} \hat{g}_\t(y) > 0$ for all $y \in [b, b + \alpha]$ and $\t_0 \in [0, 1]$. Arguing analogously for the point $a$, we will then assume that ${\left.\ddt\right|_{\t = \t_0} \hat{g}_\t(y) > 0}$ for $y \in [a - \alpha, a]$. In restriction to the arc $A = [b + \alpha, a - \alpha]$, the family~$\hat{g}_\t$ is a family of true Poincar{\'e} maps and therefore is close to the restriction to this arc of the family $\hat{f}_\t$ of Poincar{\'e} maps for the family~$V$ written in the new coordinates. The coordinate change is close to the identity and has small derivative in the parameter, so for~$\hat{f}_\t$, and hence also for~$\hat{g}_\t$, for $y \in A$ and every $\t$ the derivative in the parameter is positive.

The only thing left is to prove convergence~\eqref{eq:conv}. A chart on a transversal to a local stable or unstable separatrix of a saddle is called \emph{natural} if its origin is at the point of intersection with the separatrix. If a family of vector fields is considered, this should hold for all parameter values.
Fix some parameter value~$\t_0$. Note that for the stable separatrices of $P_\theta$ natural charts can be obtained from our new coordinate chart on~$\Sigma$ by parameter-independent shifts. Further note that, obviously, a natural chart on~$\Sigma$ for the free unstable separatrix can be obtained by a parameter-dependent shift. Take these natural charts in the upper semi-neighborhoods of~$b$ and $c_\t$. Denote by~$\lambda(\t)$ the characteristic value of the saddle~$P_\t$, i.e., the absolute value of the ratio between the negative eigenvalue and the positive one, and denote by~$\Delta_\t$ the monodromy map from the upper semi-neighborhood of~$b$ to the upper semi-neighborhood of~$c_\t$ in the natural charts. According to~\cite[Lemma~5]{IKS}\footnote{This lemma requires the natural charts to be $C^3$ both in~$y$ and the parameter. For this reason we insisted that our coordinate change was~$C^3$-smooth.}, there is a neighborhood of~$\t_0$ where for the monodromy map $\Delta_\t$ we have

%\be\label{eq:IKS1}\Delta_\t(y) = O(y^{\lambda(\t)}),\ee
%\be\label{eq:IKS2}(\Delta_\t)'_y(y) = O(y^{\lambda(\t) - 1}),\ee
\be\label{eq:IKS3}(\Delta_\t)'_\t(z) = O(z^{\lambda(\t)}\log z), \text{ as } z\to +0,\ee
where the constant that is implicitly present in the $O$-notation is independent of~$\t$.

In a semi-neighborhood $U_b$ of the point~$b$ the map~$\hat{g}_\t$ can be written as
\be\label{eq:compos}
\hat{g}_\t(y) = \Delta_\t(y - b) + \hat{c}(\t).
\ee
By \eqref{eq:IKS3} and~\eqref{eq:compos}, we have, as $y\to b + 0$,
$$\ddt\hat{g}_\t(y) = \hat{c}'(\t) + (\Delta_\t)'_\t(y - b) = \hat{c}'(\t) + o(1),$$
where the convergence in $o(1)$ is uniform in~$\t$ in the aforementioned neighborhood of~$\t_0$. This gives us the required locally uniform convergence for the derivative.
For the point~$a$ the argument is analogous. The proof of Proposition~\ref{prop:monot} is complete.

\section{Proof of Theorem~\ref{thm:main}}

\subsection{Theorems of Boyd and Veerman}

We prove Theorem~\ref{thm:main} with the help of the following results of C.~Boyd and J.~P.~P.~Veerman.

\begin{thm}[Veerman,~{\cite[Theorem~6.4]{V}}]\label{thm:veer}
Suppose that the family $\{f_\t\}_{\t \in I}$ of circle maps satisfies the following conditions:
\begin{enumerate}[label={[\arabic*]}]
\item for every $\t$ the map $f_\t$ has degree one and preserves the orientation;
\item $f_\t(x) \in C^0$ as a function of $(\t,x)$;
\item\label{cond:3} for every $\t$ the map $f_\t$ is constant on some segment~$U_\t$ that contains some point $O$ that does not depend on the parameter;
\item\label{cond:4} for any $\t_0 \in I$, for every $x_0 \in S^1\setminus \partial U_{\t_0}$ there exists a neighborhood of~$\t_0$ where the map $\t \mapsto f_\t(x_0)$ is $C^1$-smooth and the map $\t \mapsto \frac{\partial}{\partial x}f_\t(x_0)$ is~$C^0$;
\item for any point $x_0$, $f_\t(x_0)$ is non-decreasing in~$\t$, and if $x_0 \in \mathrm{int}(U_{\t_0})$ one has $\left.\frac{\partial}{\partial \t}\right|_{\t = \t_0}f_\t(x_0) > 0$.
\item\label{cond:exp} for any $\t$ the map $f_\t$ expands uniformly outside the segment $U_\t$, i.e.,
$$\forall x_0\notin U_\t, \;\;\;\; \left.\frac{\partial}{\partial x}\right|_{x= x_0}f_\t(x) > \gamma > 1,$$
where the constant $\gamma$ is independent of~$\t$.  
\end{enumerate}
Then $\mathrm{dim}_H(\{\t\colon \rho(f_\t)\notin\QQ/\ZZ\}) = 0$. Here $\mathrm{dim}_H$ is Hausdorff dimension. 
\end{thm}

\begin{rem}
This statement slightly differs from the one in~\cite{V}. Section 2 of~\cite{V} has no condition~\ref{cond:exp}, but has the requirement that for all values of the parameter the map $f_\t$ can be extended as a local diffeomorphism to a neighborhood of the arc $S^1 \setminus U_\t$, the function $\log\frac{\partial}{\partial x}f_\t(x)$ being of bounded variation. This latter requirement is, of course, not satisfied for families that we consider, but this is not a problem. In~\ref{cond:exp} this condition is utilized to prove via a Denjoy-like argument that the maps of the family have no wandering intervals outside the preimages of the flat segment. In the version above this immediately follows from condition~\ref{cond:exp}, because such an interval would grow exponentially when we iterate the map. The proof of this version of the theorem is contained in sections~5,~6 of~\cite{V}. Note that expansiveness (rather, a slightly weaker property of $(\gamma, m)$-expansiveness) is used in the proof substantially, whereas boundedness of the derivative in $x$ is not used at all. Moreover, as Veerman himself remarks, this version of the theorem is essentially due to Boyd~\cite{B}: Boyd's paper contains the proof for the case of his special families, but it can be generalized to the assumptions above.
\end{rem}

The following result is proven in~\cite{B} at pp.~44--46 as part of the proof of Boyd's Theorem~3.

\begin{thm}[C. Boyd, {\cite[Section~5]{B}}]\label{thm:superboyd}
Let $V$ and $W$ be two families of vector fields of class~$\N$ and $\{f_\t\}_{\t\in[0,1]}$, $\{g_\tau\}_{\tau\in[0,1]}$ be the corresponding families of Poincar{\'e} maps. Suppose there exists a homeomorphism $s \colon [0, 1] \to [0, 1]$ such that for any $\t \in [0, 1]$ one has
$\rho(g_{s(\t)}) = \rho(f_\t)$.
Assume also that in both families the rotation number depends monotonically on the parameter, there are no (non-degenerate) segments in the parameter space that correspond to an irrational rotation number, and that for every parameter segment where the rotation number is rational there are no interior points that correspond to fields with saddle loops\footnote{In particular, there must be no segments that correspond to fields with separatrix loops.}. Then the families~$V$ and~$W$ are strongly equivalent.
\end{thm}
\begin{rem}
In this theorem, the condition of monotonic dependence of the rotation number on the parameter can be omitted. Instead of requiring that there are no segments with irrational rotation number or interior points that correspond to the presence of loops one can require that $s$ establish a correspondence between such segments and points for the two families.
\end{rem}

\subsection{The main proposition}
Roughly speaking, the following proposition claims that for a deformation of any vector field $v\in \orc(M)$ one can choose a neighborhood $U$ (which is a torus without a disk, i.e., a handle) of the unstable manifold~$W^u(P_v)$ in such a way that the restriction of the deformation onto this neighborhood can be extended as a family of vector fields on the whole torus so that after reversing the time one gets a family of fields of class~$\N$ with some additional properties that will turn useful in proving Theorem~\ref{thm:main}. The most important property is the monotonicity in the parameter, because it will allow us to apply Theorems~\ref{thm:veer} and~\ref{thm:superboyd}.

\begin{prop}\label{prop:immer}
Let $M$ be an arbitrary surface with a handle, $v\in\orc(M)$, and $V = \{v_\t\}_{\t\in[-1, 1]}$ be a generic family of smooth vector fields on~$M$ such that $v_0 = v$.
Then there exists a neighborhood $U \subset M$ of the unstable manifold $W^u(P_v)$, homeomorphic to a torus without a disk, a segment $J\ni 0$ in the parameter space, and a smooth embedding $\imath\colon U \to \TT$ such that the family of vector fields $\{\imath_*(-v_{\t})\}_{\t \in J}$ can be extended to the whole torus~$\TT$ as a family ${\tilde{V} = \{\tilde{v}_\t\}_{\t \in J}}$ with the following properties.
\begin{enumerate}
\item The fields of the family ${\tilde{V} = \{\tilde{v}_\t\}_{\t \in J}}$ are of class~$\N$.

\item For the corresponding family $F = \{f_\t\}$ of Poincar{\'e} maps the following holds:

  \begin{enumerate}\setlength\itemsep{-0.1em}
    \item $\rho(f_0) = 0$;
    \item the function $r(\t) = \rho(f_{\t})$ assumes rational values on some non-degenerate (i.e., not equal to a point) segments; the fields of the family have saddle loops exactly for those parameter values that are the endpoints of these segments;
    \item there exists a finitely-smooth parameter-depending coordinate change such that after this change and, maybe, after changing the sign of the parameter one would have, for all~$\t$ and~$y$, the inequality $\ddt \hat{f}_{\t}(y) > 0$, where $\hat{f}_{\t}$ is $f_{\t}$ written in the new coordinate.
  \end{enumerate}

\item $\Bif(\tilde{V}) = \overline{\{\t\colon \rho(f_\t)\notin\QQ/\ZZ\}} = \overline{\{\t\colon \tilde{v}_\t  \text{ has a loop}\}}$, and $\Bif(\tilde{V})$ is a Cantor set.
\end{enumerate}

\begin{rem}\label{rem:close}
1. It will be clear from the proof that for a family~$V_1$ sufficiently close to~$V$ we can take the same sets~$U, J$ and embedding~$\imath$.

2. The neighborhood $U$ will be trapping for $v_\theta$, so the transitive quasi-minimal sets that appear in $\bar{V}$ for irrational rotation numbers will be contained in $\imath(U)$ and so such sets also appear in the original family~$V$ at a set of parameters of continuum cardinality.
\end{rem}
\end{prop}
We postpone the proof of this proposition until section~\ref{sect:immer} and first deduce Theorem~\ref{thm:main} from the proposition.

\subsection{A restatement of theorem~\ref{thm:main}}

Before we start proving Theorem~\ref{thm:main}, let us restate it in a more precise way and at the same time introduce some notation that will be used in the proof.

\begin{mytheorem}\label{thm:main2}
Let $v\in\orc(M_1)$ and $w\in\orc(M_2)$. For generic smooth families $V = \{v_\t\}_{\t \in [-1, 1]}, \; v_0 = v,$ and $W = \{w_\t\}_{\t \in [-1, 1]}, \; w_0 = w,$ there exist segments $J_1, J_2$: $0\in J_j \subset [-1, 1],$ and neighborhoods $U_1\subset M_1, \, U_2 \subset M_2$ of $W^u(P_{v_0})$ and $W^u(P_{w_0})$ respectively such that for the restrictions $\overline{V} = V|_{J_1 \times U_1}$ and $\overline{W} = W|_{J_2 \times U_2}$ we have the following.
\begin{itemize}
\item The families $\overline{V}$ and $\overline{W}$ are strongly equivalent; each is strongly structurally stable.

\item The bifurcation diagram $\Bif(\overline{V})$ is a Cantor set~$K$ that contains 0 and lies on one side from it. The boundaries of the intervals of $J_1\setminus K$ correspond to vector fields with a separatrix loop for~$P$. Whenever there is a separatrix loop for~$P$, the free unstable separatrix winds onto it.

\item $\Bif(\overline{V})$ has Hausdorff dimension zero (and, therefore, zero Lebesgue measure).

\end{itemize}
\end{mytheorem}

\subsection{Equivalence for cropped families}

The genericity condition implicitly imposed on families in both Theorem~\ref{thm:main2} and Proposition~\ref{prop:immer} is that the saddle loop must be unfolded with non-zero speed as the parameter changes. This implies that on one side of the zero value a repelling hyperbolic cycle is born from the loop. 
Applying Proposition~\ref{prop:immer} to the families $V$ and $W$ from the statement of Theorem~\ref{thm:main2}, we get two families $\tilde{V}$ and $\tilde{W}$ of vector fields on the torus (and also corresponding segments $J_1, J_2$ and embeddings $\imath_1,~\imath_2$).

We want to apply Theorem~\ref{thm:superboyd} to these families $\tilde{V}$ and~$\tilde{W}$. First, note that by Proposition~\ref{prop:immer} the Poincar{\'e} maps for each family become strictly monotonic in~$\t$ when written in some parameter-dependent chart. Since the rotation number does not depend on the chart, this implies that the rotation number depends non-strictly monotonically on~$\t$ in general and is strictly monotone at points where irrational values are assumed~\cite[Lemmas~1-3]{Ha} (see also \cite{He} and \cite[Prop, 11.1.8-9]{KH}), so there are no segments where the rotation number is irrational. By the proposition, saddle loops are present only for parameter values at the endpoints of segments where the rotation number is rational. Since the zero value of the parameter corresponds to zero rotation number and to the presence of a loop, we can change the sign of the parameter if necessary and assume that in the left semi-neighborhood of zero the rotation number equals zero. Furthermore, after applying, if necessary, a diffeomorphism of the torus that changes orientation on the transverse circle~$\Sigma$, we can assume that the rotation number does not decrease with the parameter. We can now crop the parameter segments $J_1$ and $J_2$ in such a way that after that for the families $F$ and $G$ of Poincar{\'e} maps the rotation number will run through the same segment and the endpoints of the cropped segments will correspond to structurally stable vector fields.
Now we can construct a homeomorphism between the cropped~$J_1$ and~$J_2$ as follows: for each pair of segments where the rotation number is the same,
let~$s$ take the first segment onto the second one affinely. Since the rotation number for both families depends on the parameter continuously and monotonically, the map~$s$ can be extended as a homeomorphism between the parameter spaces.

Theorem~\ref{thm:superboyd} yields the strong equivalence of the families $\tilde{V}$ and~$\tilde{W}$. Let $H\colon (\t, z) \mapsto (h(\t), H_\t(z))$, $z\in\TT$, be a homeomorphism that realizes this equivalence. Now we can establish the equivalence of the restrictions from the first statement of Theorem~\ref{thm:main2}. Take a neighborhood $T_1 \subset \imath_1(U_1)$ of the unstable manifold $W^u(P_{\tilde{v}_0})$ such that for $\t$ near zero (i.e., for~$\t$ from some segment $\tilde{J}_1\subset J_1$ that contains zero) we have
$H_\t(T_1)\subset \imath_2(U_2).$
This is almost what is required, except in the image the neighborhood of the unstable manifold depends on the parameter.

Assume that for the cropped parameter segment~$\tilde{J}_1$ the endpoints still correspond to structurally stable fields.
Denote $\imath_2(U_2)$ by~$T_2$. We want to modify the restriction $H|_{\tilde{J}_1 \times T_1}$ in the neighborhood of the set $\tilde{J}_1 \times \partial T_1$ in order to make the modified map a homeomorphism onto the image that, for each $\t\in\tilde{J}_1$, realizes the equivalence between $\tilde{v}_\t|_{T_1}$ and~$\tilde{w}_\t|_{T_2}$.

Let $A$ be an annular neighborhood of $\partial T_2$ in~$T_2$ such that $H_\t(\partial T_1) \subset A$ for all $\t \in \tilde{J}_1$. We can simultaneously rectify all the restrictions~$\tilde{w}_\t|_A$ by a change of coordinates that is continuous and even smooth in the parameter. Thus, we can think of~$A$ as the annulus $\{(\varphi, r) \mid 1 \le r \le 2\}$ in polar coordinates, with trajectories being radial segments pointing towards zero. The curves $H_\t(\partial T_1)$ are $\t$-dependent topological circles which together with the circle~$\{r = 1\}$ bound smaller topological annuli. Clearly, for any~$\t$ there is a homeomorphism $G_\t$ of the form
$(\varphi, r) \mapsto (\varphi, f_\t(r))$
that takes the annulus in-between the circles~$\{r = 1\}$ and~$H_\t(\partial T_1)$ to~$A$ and equals the identity at~$\{r = 1\}$. Moreover, we can make homeomorphisms $G_\t$ depend continuously on~$\t$. Then the required modification of $H|_{\tilde{J}_1 \times T_1}$ can be made by post-composing each $H_\t$ with~$G_\t$ on $H_\t(T_1) \cap A$.

Finally, we have established the first claim of the theorem, with $\tilde{J}_1$ playing the role of~$J_1$, $h(\tilde{J}_1)$ being~$J_2$, and the neighborhoods of the unstable manifolds being equal to $\tilde{U}_1 = \imath_1^{-1}(T_1)$ and $\tilde{U}_2 = \imath_2^{-1}(T_2)$, respectively.

The strong structural stability of the cropped families is proved analogously using Remark~\ref{rem:close}.

\subsection{The bifurcation diagram}
By Proposition~\ref{prop:immer}, for the family $\tilde{V}$ the bifurcation diagram coincides with the set~$\overline{\{\t\colon \tilde{v}_\t  \text{ has a loop}\}}$ and is a Cantor set.
Moreover, Theorem~\ref{thm:veer} can be applied to the family~$\{\hat{f}_\t\}$ of its Poincar{\'e} maps written in certain parameter-dependent coordinates.
As we will see in the proof of Proposition~\ref{prop:immer}, the flat interval of the map $\hat{f}_\t$ does not depend on the parameter, but were it not the case, we could argue that it depends continuously on the parameter and so condition~\ref{cond:3} of Theorem~\ref{thm:veer} is satisfied after further cropping the parameter space.

Theorem~\ref{thm:veer} yields that $\dim_H(\Bif(\tilde{V}|_{\tilde{J}_1})) = 0$, if the segment $\tilde{J}_1$ is sufficiently small. For the family $V|_{\tilde{J}_1, \tilde{U}_1}$, the bifurcation diagram is the same.

\subsection{The free separatrix winds onto the loop}
For a vector field of class~$\N$ that has a saddle loop, the free separatrix of the saddle~$P$ has to wind onto this loop. Indeed, on the one hand, in reversed time the loop attracts every point in its narrow monodromic semi-neighborhood. On the other hand, consider a small segment $L \subset \Si$ that has its endpoint at the loop and is contained in this semi-neighborhood. Since the Poincar{\'e} map $f$ is expansive outside the flat interval, the images~$f^k(L)$ have to grow exponentially until one of them intersects the flat interval and hence intersects the free separatrix of the saddle. This implies that the free separatrix is attracted to the loop in reversed time, that is, it winds onto the loop. Proposition~\ref{prop:immer} yields that the same holds for the loops of the saddle~$P$ in family~$V$.
Theorem~\ref{thm:main2} (or~\ref{thm:main}) is proven modulo the proof of Proposition~\ref{prop:immer}.

\subsection{Proof of Proposition~\ref{prop:immer}}\label{sect:immer}

\subsubsection*{Idea of the proof}
The set $U$ from the statement of Proposition~\ref{prop:immer} is a narrow trapping neighborhood of~$W^u(P)$ for the field $v = v_0$. It is a torus-without-a-disk embedded into the surface~$M$. We take $U$ with the field~$v$ attached to it and glue to the boundary of $U$ a disk that supports a vector field with a single repeller. This gluing gives us an embedding~$\imath_U$ of the surface~$U$ into a torus~$\TT$. After that we invert the direction of the vector field $(\imath_U)_*(v)$ and start modifying  it outside the image~$\imath_U(U)$ in order to obtain a vector field of class~$\N$. Then we pushforward the rest of the family~$V$ and obtain, after the same procedure and appropriate cropping, a small family of fields of class $\N$ on the whole torus. Then we modify this family to make the corresponding Poincar{\'e} map depend monotonically on the parameter. This latter property yields the desired facts about the bifurcation diagram.

\subsubsection*{Constructing $\tilde{v}_0$ of class~$\N$}
First, take a narrow trapping neighborhood~$U$ of the unstable manifold~$W^u(P)$ for the field~$v$ (we will specify later how narrow it should be). This neighborhood must be a torus without a disk, the boundary circle~$\partial{U}$ being transverse to the vector field. Take also a unit disk~$D$ with a field $x\frac{\partial}{\partial x} + y\frac{\partial}{\partial y}$. Note that this field has a unique singularity which is a hyperbolic source and this field is transverse to the boundary. Glue~$U$ to~$D$ along the neighborhoods of the boundaries in such a way that the vector fields be glued as well to form a smooth vector field on the torus~$\TT = U \# D$. Multiply this field by~$-1$ and denote the resulting field on the torus by~$\tilde{w}$. Multiplying by~$-1$ is equivalent to inverting the time. So, the repeller on the glued-in disk becomes an attractor, which we will denote by~$\Omega$. The separatrices of the saddle~$P$ change their directions. The embedding $\imath_U\colon U \to \TT$ that originates from gluing is exactly the one that appears in the statement of Proposition~~\ref{prop:immer}, and~$U$ is the same, provided it was chosen sufficiently narrow. Denote by $\imath_D$ the embedding of~$D$ in~$\TT$; denote $\imath_U(U)$ by~$\tilde{U}$ and~$\imath_D(D)$ by~$\tilde{D}$.

The field~$\tilde{w}$ has a global closed non-contractible transverse circle. To obtain one, first take a point at the local unstable separatrix of the saddle~$P$ that is involved into the loop and draw a transverse curve through it so that both endpoints of this curve are on the boundary of~$\tilde{U}$. Denote these two points by~$p_1, p_2$. Since the original field on~$D$ was radial, it is clear that~$p_1$ and~$p_2$ can be connected by a transverse curve that goes inside~$\tilde{D}$ so that this curve together with the previous curve between~$p_1$ and~$p_2$ makes a smooth global transversal. However, one can obtain in this manner transversals which are of different homological type. In order to make the relative arrangement of~$\tilde{U}$ and the transversal as simple as possible, we do the following. First,~$\tilde{U}$ is cut into three parts: a neighborhood of the saddle and two strips attached to it; the strips may be viewed as neighborhoods of a segment of the loop and a segment of the free separatrix. Then the transverse segment $p_1p_2$ is taken inside the strip that is a neighborhood of the segment of the loop. Denote the closed transversal obtained from this segment by~$\Si$.

Now take a smooth chart (that ``unfolds'' the torus into a rectangle) in which the circle~$\Si$ becomes a vertical segment and the field in the neighborhood of~$\Si$ becomes horizontal. Inside this chart, draw another vertical transverse circle~$\hSi$ near~$\Si$ (see Fig.~\ref{fig:v0}) so that the unstable separatrix involved into the loop cross it before it crosses~$\Si$. In what follows we will refer to the geometry of the picture when we say ``above'' or ``below'', etc.
Denote by $O$ and $\hat{O}$ the points where $W^u_{loc}(P)$ first intersects~$\Si$ and~$\hSi$ respectively.

\begin{figure}
  \begin{center}

  \begin{tikzpicture}[scale=2]

  % domain where we can perturb the field, filled light green 
  \fill [fill = green!20!white] (2.8,4) -- (2.8, 2.6) .. controls (3.25,2.49) .. (3.9,2.45) -- (3.9, 4);
  \fill [fill = green!20!white] (2.8,0) -- (2.8, 1.4) .. controls (3.24,1.515) .. (3.9,1.55) -- (3.9, 0);

  % points d_1, d_2
  \filldraw [color = black] (2.8, 2.6) circle (0.5pt) node[below right] {{\scriptsize $d_2$}};
  \filldraw [color = black] (2.8, 1.4) circle (0.5pt) node[above right] {{\scriptsize $d_1$}};

  %separatrices
  \coordinate (P) at (1.5, 2.2);
  \draw[violet, very thick] plot [smooth, tension = 0.7] coordinates{(0,2) (1.1,2) (P)};
  \draw[red, very thick]  plot [smooth, tension = 0.7] coordinates{(P) (1.7, 2.8) (1.8,4)};
  \draw[violet, very thick] plot [smooth, tension = 0.7] coordinates{(P) (2.3,2) (5,2)};
  \draw[blue, very thick] plot [smooth, tension = 0.9] coordinates{(P) (0.7, 3) (0.05,3.95)};
  \draw[red, very thick]  (1.8, 0) .. controls (1.8, 1.6) and (0.9, 1.65) .. (0,1.75);
  \draw[red, very thick, dash pattern=on 2 off 2 on 4 off 2 on 5 off 2 on 1000]  (1.55, 1.81) .. controls (3, 1.77) and (4, 1.78) .. (5,1.75);

  %arrows on separatrices
  \draw[->, violet, thick] (0.8, 1.985) -- (0.9, 1.99);
  \draw[->, violet, thick] (2.2, 2.01) -- (2.3, 2);
  \draw[->, blue, thick] (0.8, 2.88) -- (0.7, 3.);
  \draw[<-, red, thick] (1.65, 2.58) -- (1.675, 2.68);
  \draw[->, red, thick] (2.2, 1.794) -- (2.3, 1.79);
  \draw[<-, red, thick] (1.4, 1.3) -- (1.3, 1.4);

  %saddle P
  \filldraw [color = black] (P) circle (1pt) node[below] {$P$};
  
  %transverse circles
  \draw[thick] (2.8,0) -- (2.8,4) node[above, color = black] {$\hat\Sigma$};
  \filldraw [color = black] (2.8, 1.985) circle (0.5pt) node[above right] {{\scriptsize $\hat{O}$}};
  \draw[thick, color = black] (3.9,0) -- (3.9,4) node[above, color = black] {$\Sigma$};
  \filldraw [color = black] (3.9, 1.985) circle (0.5pt) node[above right] {{\scriptsize $O$}};

  %neighborhood $\tilde{U}$
  \node at (2.1, 3.8) {$\tilde{U}$};
  \draw (2.3,4) .. controls (2.3,2.5) .. (5,2.4);
  \draw (2.3,0) .. controls (2.3,1.5) .. (5,1.6);
  \draw (1.2,4) .. controls (1.2,2.7) and (1.1, 2.45) .. (0,2.4);
  \draw (1.2,0) .. controls (1.25,1.45) and (0.8, 1.54) .. (0,1.6);

  %points q, \hat{q}
  \filldraw [color = black] (2.8, 1.78) circle (0.5pt) node[below left] {{\scriptsize $\hat{q}$}};
  \filldraw [color = black] (3.9, 1.77) circle (0.5pt) node[below left] {{\scriptsize $q$}}; 

  %torus label
  \node[above] at (0, 4) {$\mathbb{T}^2$};

  %Omega
  \filldraw [fill = red, draw = red] (0,3.92) arc (-90 : 0 : 0.08) -- (0,4);
  \node at (0.23, 3.9) {$\Omega$};
  \filldraw [fill = red, draw = red] (0,.08) arc (90 : 0 : 0.08) -- (0,0) node[above right] {$\Omega$};
  \filldraw [fill = red, draw = red] (4.92,4) arc (180 : 270 : 0.08) -- (5,4) node[below left] {$\Omega$};
  \filldraw [fill = red, draw = red] (4.92,0) arc (180 : 90 : 0.08) -- (5,0) node[above left] {$\Omega$};

  % transverse circle around Omega
  %\draw (coordinate) arc (start angle : end angle : radius)
  %\draw[red, dotted] (0.5,0) arc ( 0 : 90 : 0.5);
  %\draw[red, dotted] (4.5,0) arc ( 180 : 90 : 0.5);
  %\draw[red, dotted] (0.5,4) arc ( 0 : -90 : 0.5);
  %\draw[red, dotted] (4.5,4) arc ( 180 : 270 : 0.5);

  %torus rectangle
  \draw[thick] (0,0) rectangle (5,4);
  \end{tikzpicture}
  \caption{The structure of the fields~$\tilde{w}$ and $\tilde{v}_0$. The domain where we can alter the field is shaded.}\label{fig:v0}
  \end{center}
\end{figure}

In the lower semi-neighborhood of the point~$O$ on~$\Sigma$ the monodromy map $\Delta\colon\Si \to \hSi$ is defined. It is expanding in restriction to a sufficiently small lower semi-neighborhood~$B \subset \Si$ of~$O$, because the saddle is area-expanding now. Our monodromy map can be continued to the point~$O$ by specifying that~$O$ be mapped into~$\hat{O}$. Furthermore, $O$ is accumulated from below by a sequence of intersection points of~$\Si$ and the free stable separatrix of the saddle. Take the very first point of this sequence and denote it by~$q$; denote by~$\hat{q}$ the analogous point on the second transverse circle. We want the $\Delta$-image of~$B$ to contain~$\hat{q}$. It can be achieved by choosing the neighborhood~$U$ sufficiently narrow in the first place, because for the field~$v$ the monodromy map along the saddle loop was a strong contraction, provided that we looked at a sufficiently narrow monodromic semi-neighborhood of the loop. Thus, we can and will assume that the $\Delta$-image of $B$ has an endpoint $d_1$ at the boundary of~$\tilde{U}$.

In the upper semi-neighborhood $\Gamma\subset\Si$ of the point~$q$ the monodromy map to~$\hSi$ is defined too. Note, first, that it can be continued to~$q$ by letting $q$ be taken to~$\hat{O}$. The monodromy map takes~$\Gamma$ into an upper semi-neighborhood of the point~$\hat{O}$. Moreover, if~$\Gamma$ is sufficiently small, this map is expanding due to crossing the hyperbolic sector of the saddle~$P$, where the local monodromy map is expanding since the saddle is area-expansive. We will assume that we have chosen~$U$ so narrow that the image of~$\Gamma$ under the monodromy map has an endpoint $d_2\in\partial(\tilde{U})$ and the monodromy map is expansive on~$\Gamma$. We can do that because after choosing the segment~$\Gamma$ where the monodromy map is expansive, we could look at its monodromy image and crop our neighborhood~$\tilde{U}$ to make the endpoint of the image be at the boundary of the cropped neighborhood. The set $U$ in the preimage is also cropped. After this cropping, the monodromic semi-neighborhood of the loop contained in~$U$ does not change, only the neighborhood of the free local separatrix and the non-monodromic semi-neighborhood of the loop become narrower, so we can assume that we have chosen~$U$ to be like that.

Thus, we have an expansive monodromy $\Delta \colon B\cup\Gamma \to [d_1, d_2]\subset~\hSi$. Note that on the arc $(O, q)$ (recall that according to our notation this arc goes from~$O$ upwards and then approaches~$q$ from below) the monodromy map is not defined because every orbit that starts at this arc goes to the sink~$\Omega$. So, we extend the map to $(O, q)$ by setting $\Delta((O, q)) = \{\hat{O}\}$, as usual. Further note that $\Delta$ continues to the arc $(q, O)$, but it does not have to be expansive throughout this arc. We want to obtain a field that belongs to the class~$\N$, and the only thing we need for that is to make the Poincar{\'e} map expansive outside its flat segment. To achieve that, we modify the field in between~$\hSi$ and~$\Si$ using a method from~\cite[Lemma~2, p.~183]{PM}. That is, we consider an expansive smooth map $f\colon [q, O] \to \hSi$ that, when written in our coordinates, coincides with the monodromy map $\Delta$ on $B\cup\Gamma$. It is clear that such~$f$ exists.
We denote the coordinate representations of our maps by the same symbols that we use for the maps themselves. We set
$$\varphi(y) = f \circ \Delta^{-1}(y),$$
$$\varphi_s(y) = (1 - s)y + s\varphi(y), \;\; \forall s \in [0, 1],$$
$$H(x, y) = (x, \varphi_{\sigma(x)}(y)),$$
where $\sigma$ is a smooth monotonic function that takes values in~$[0, 1]$, is equal to zero in a neighborhood of the $x$-coordinate~$x_\hSi$ of~$\hSi$, and is equal to one in a neighborhood of the $x$-coordinate~$x_\Si$ of~$\Si$.

Finally, we set
$$\tilde{v}_0(x,y) = DH(H^{-1}(x,y))[\tilde{w}(H^{-1}(x,y))]$$
(the derivative of $H$ taken at the point $H^{-1}(x,y)$ is applied to the vector of the field $\tilde{w}$ taken at the same point) in the strip between~$\hat\Sigma$ and~$\Sigma$ and set~$\tilde{v}_0 = \tilde{w}$ outside this strip. It is easy to check that for the field~$\tilde{v}_0$ the monodromy from~$\hSi$ to~$\Si$ coincides in coordinates with the quotient~$\varphi = f \circ \Delta^{-1}$, and so the (extended) Poincar{\'e} map from~$\Si$ to itself is well defined and coincides in coordinate representation with the expanding map~$f$ on~$[q, O]$. Hence we have obtained the Poincar{\'e} map that expands everywhere outside the segment where it is constant. It is also easy to check that we did not alter the field inside~$\imath_U(U)$: the maps~$f$ and~$\Delta$ coincide on  $B\cup\Gamma$, therefore $\varphi_s(y) = y$ for all~$s \in [0, 1]$ and $y\in\Delta(B\cup\Gamma)$, which implies that for the points~$(x,y)$ in the intersection of $\tilde{U}$ and the vertical cylinder $C$ between $\hat\Si$ and $\Si$ we have $H^{-1}(x,y) = (x, y)$, and hence $DH = {Id}$, so the field remains horizontal there.

It is now clear that the vector field $\tilde{v}_0$ is of class~$\N$.

\subsubsection*{Monotonicity in the parameter}
Consider once again the original family $V$ on~$M$. When the parameter is close to zero, the field~$v_{\t}$ differs from~$v_0$ only very slightly in a neighborhood of~$\partial U$ and, therefore, is transverse to~$\partial U$, so we can extend the fields $(\imath_U)_*(-v_\t)$ simultaneously to obtain a smooth family $W = \{w_\t\}_{\t\in[-\e,\e]}$ of vector fields defined on the whole torus, with $w_0~=~\tilde{v}_0$, by gluing, as above, a disk with a single repeller, but now by a para\-meter-depen\-dent correspondence between the neighborhoods of the boundaries. No matter how we choose this correspondence, if we take~$\e$ sufficiently small, we will have~$w_\t\in\N$ for all~$\t$, because $w_0 = \tilde{v}_0 \in \N$ and~$\N$ is open.

Let $G$ be the family of the Poincar{\'e} maps for the family~$W$, where the transverse circle is still~$\Si$.
Our goal now is to modify the family~$W$ so that the modified Poincar{\'e} maps become monotone in the parameter, at least in some parameter-dependent coordinate. 

In a narrow neighborhood of~$\Si$ the vector field can be rectified, simultaneously for all parameter values, into a unit horizontal field by a parameter-dependent smooth change of coordinates that leaves~$\Si$ vertical. We switch to this parameter-dependent chart and draw another vertical transversal~$\hSi$ near $\Si$. Since the chart depends on the parameter, so does~$\hSi$ and also the strip by which $\tilde{U}$ intersects the cylinder~$C$ which is in between the transverse circles. Define $O, \hat{O}, q, \hat{q}, B, \Gamma, d_1, d_2$ as above. From now on, the coordinates on~$\Si$ and~$\hSi$ will only be changed in accord, so that we could always assume that the vector field in~$C$ is unit horizontal.

It is time to refer to the genericity condition. Consider a parameter-dependent coordinate on the transversal~$\Si$ such that its origin for every parameter value coincides with the point of first intersection between $\Si$ and the stable separatrix involved in the loop when the parameter is zero. Recall that such coordinates are called natural with respect to this separatrix. The genericity condition is that the point~$c(\t)$ of the first intersection of~$\Si$ and the unstable separatrix involved in all the loops has non-zero derivative in the parameter. Note that, when we switch to a different natural chart, the sign of the derivative does not change if two charts are co-oriented. After changing, if necessary, the sign of the parameter~$\t$, we will assume that~$c'(0) > 0$ and, therefore, for~$\t$ near zero one also has~$c'(\t) > 0$.

We would like to work in a natural chart (with respect to the same stable separatrix) such that in it the point of the first intersection between~$\Si$ and the other stable separatrix does not depend on the parameter. When we switch to this chart, inequality $c'(\t) > 0$ still holds for~$\t$ near zero.

Now, let us assume that we have chosen the neighborhood~$U$ so narrow that in our natural chart we now have

$$\left.\ddt\right|_{\t = 0}g_\t(y) > \delta > 0, \;\;\; \forall y\in B\cup\Gamma.$$

It is indeed possible to have done that. Arguing as in the proof of Proposition~\ref{prop:monot}, we can conclude that in the upper semi-neighborhood of the point~$q$ and in the lower semi-neighborhood of~$O$ the derivative under consideration is positive. Concerning the lower semi-neighborhood of the point~$O$, when choosing~$U$, we could ensure that the image of this semi-neighborhood under the monodromy map from~$\Si$ to itself have one endpoint at~$\partial \tilde{U}$. For the segment~$\Gamma$ we can argue analogously: take the upper semi-neighborhood of~$q$ where the derivative in the parameter is positive, consider the image of the upper endpoint of this semi-neighborhood under the monodromy map, and crop~$\tilde{U}$ in such a way that this point is at the new boundary.

We can now modify the family. Consider the cylinder~$C$ and the coordinates~$x,y$ on it in which all vector fields of our family appear constant and horizontal. Let 

\begin{itemize}

\item $\sigma(x)$ be a $C^\infty$-smooth bump-function that has support contained in $[x_\hSi, \, x_\Si]$, takes values in $[0, \, 1]$, and is equal to 1 at the point~$(x_\hSi + x_\Si)/2$;
\item $\kappa(y)$ be a $C^\infty$-smooth bump-function equal to 1 exactly in the $\frac{\alpha}{2}$-neighborhood of the segment $[d_1, d_2]$ and equal to zero outside its $\alpha$-neighborhood, where $\alpha > 0$ is a small constant that we will specify below.

\end{itemize}
Let us add to our family the vector field
$$\Phi = K\t\cdot\sigma(x)(1 - \kappa(y))\cdot\frac{\partial}{\partial y} + 0\cdot\frac{\partial}{\partial x} + 0\cdot\frac{\partial}{\partial \t},$$
where $K > 0$ is a large constant, and see what effect it has on the derivative of the Poincar{\'e} map in the parameter.

Denote by $h_\t(y)$ the parameter-dependent monodromy map from $\hSi$ to~$\Si$ for the modified family, written in the same coordinates on~$C$. For notational convenience, we will also write $h(y, \t) = h_\t(y)$, and likewise for all parameter-dependent functions. Due to the choice of the chart, we have $h_0 = Id$. Moreover, for~$\t$ near zero the map~$h_\t$ is close to the identity map. It is easy to check that outside the $\frac{\alpha}{2}$-neighborhood of the segment $[d_1, d_2]$ the derivative $(h)'_\t(y, \t)$ is non-negative and outside the $\alpha$-neighborhood of this segment it is separated from zero by some constant that depends on~$K$.

No matter what $\alpha$ we choose, we can always crop our parameter space so that the following will hold: for any parameter value, the intersection of~$\tilde{U}$ and the cylinder~$C$ lies in the rectangle $[x_\hSi, \, x_\Si] \times [d_1 - \alpha/2,\, d_2 + \alpha/2]$. Then for $\t$ near zero adding $\Phi$ does not change the vector fields inside~$\tilde{U}$. Let us now see how the derivative of the Poincar{\'e} map in $\t$ has been affected. The new Poincar{\'e}  map $f_\t$ can be written as a composition of the old one (which in our coordinates coincides with the monodromy map from~$\Si$ to~$\hSi$) and the map~$h_\t$. By the chain rule we get
$$(f)'_\t(y, \t) = \underbrace{(h)'_\t(g_\t(y), \t)}_{\mathcal{A}} + \underbrace{(h)'_y(g_\t(y), \t)}_{\mathcal{B}} \cdot \underbrace{(g)'_\t(y,\t)}_{\mathcal{C}}.$$
Here $\mathcal{B}$ is close to one, $\mathcal{C}$ is bounded from below, $\mathcal{A}$ is always non-negative and is large if~$K$ is large and $g_\t(y) \notin [d_1 - \alpha,\, d_2 + \alpha]$. Hence, the sum is positive if $g_\t(y) \notin [d_1 - \alpha,\, d_2 + \alpha]$. From~$\mathcal{C}$ we additionally require that for~$\t$ near zero it must be positive for those~$y$ that satisfy $g_\t(y)\in[d_1 - \alpha,\, d_2 + \alpha]$ --- this can be achieved by choosing~$\alpha$ sufficiently small (for $\t = 0$ the expression $\mathcal C$ is separated from zero on the $g_0$-preimage of~$[d_1, d_2]$). Thus, by taking~$K$ large and~$\alpha$ small we assure that the $\t$-derivative of the modified Poincar{\'e} map is positive and separated from zero. As soon as we choose~$\alpha$ and~$K$, the construction of the family~$\tilde{V}$ whose existence is claimed in Proposition~\ref{prop:immer} is complete.

\subsection*{Properties of the bifurcation diagram}
Vector fields of class $\N$ with irrational rotation number are not structurally stable, therefore $\overline{\{\t\colon \rho(f_\t)\notin\QQ/\ZZ\}} \subset \Bif(\tilde{V}).$ On the other hand, if a vector field of class~$\N$ has a rational rotation number and has no saddle loop, this field is Morse-Smale and therefore is structurally stable. Thus, in order to prove the inverse inclusion, it suffices to check that every parameter value that corresponds to a field with a loop lies in the set $\overline{\{\t\colon \rho(f_\t)\notin\QQ/\ZZ\}}$. In other words, it suffices to prove that in our family the rotation number, as a function of the parameter, is non-constant in a neighborhood of any parameter value that corresponds to a saddle loop. Although we have monotonic dependence of Poincar{\'e} maps on the parameter only in some specifically chosen coordinates, we can argue exactly as in section~\ref{sect:bifdiagboyd} to establish that. Of course, we also use the fact that the coordinates of the points where the stable separatrices first intersect the transversal are independent of the parameter. This latter fact and monotonicity also yield that separatrix loops happen only for parameter values at the endpoints of segments where rotation number is rational. Thus, we have equality $$\Bif(\tilde{V}) = \overline{\{\t\colon \rho(f_\t)\notin\QQ/\ZZ\}} = \overline{\{\t\colon \tilde{v}_\t  \text{ has a loop}\}}.$$ From this we deduce that $\Bif(\tilde{V})$ is perfect and nowhere dense. Since it is also obviously closed, it is a Cantor set.
This completes the proof of Proposition~\ref{prop:immer}.

\vspace{0.3cm}

\medskip
\noindent
{\large Ivan Shilin,\\}
National Research University Higher School of Economics, Russia, Moscow,\\
E-mail: \texttt{i.s.shilin@yandex.ru}

\end{document}